\documentclass{article}
\usepackage{todonotes}

\usepackage[utf8]{inputenc}

\usepackage{xcolor}
\definecolor{darkolivegreen}{rgb}{0.33, 0.42, 0.18}
\definecolor{celestialblue}{rgb}{0.29, 0.59, 0.82}

\usepackage[colorlinks,
            linkcolor=darkolivegreen,
            citecolor=darkolivegreen,
            urlcolor=darkolivegreen]{hyperref}

\hypersetup{pdftitle={CGTJS}, pdfauthor={Felipe Galarce Marin}, 
}

\usepackage{abstract}
\usepackage{authblk}


\usepackage{fullpage}
\usepackage{amsmath, amsthm, amssymb}
\usepackage{thmtools, thm-restate}
\usepackage{dsfont}
\usepackage{mathtools}
\usepackage{graphicx}   
\usepackage{subcaption}
\usepackage{color}
\usepackage{setspace}
\usepackage{acro}
\usepackage{booktabs}
\usepackage{isomath}

\renewcommand{\vec}{\vectorsym}
\newcommand{\ten}{\tensorsym}

\theoremstyle{definition}

\numberwithin{equation}{section} 

\newcommand{\cH}{\ensuremath{\mathcal{H}}}

\newcommand{\bR}{\ensuremath{\mathbb{R}}}

\def\[{\left[}
\def\]{\right]}
\def\<{\langle}
\def\>{\rangle}
\def\({\left(}
\def\){\right)}
\def\[{\left [}
\def\]{\right]}
\def\({\left(}
\def\){\right)}

\DeclareMathOperator{\grad}{\nabla}
\DeclareMathOperator{\dive}{\text{div}}
\DeclareMathOperator{\tr}{\text{tr}}

\newcommand{\bu}{\vec u}
\newcommand{\bv}{\vec v}
\newcommand{\bz}{\vec z}
\newcommand{\bw}{\vec w}
\newcommand{\bff}{\vec f}

\newcommand{\ts}{\ten \sigma}

\renewcommand{\ker}{\ensuremath{\operatorname{ker}}}
\newcommand{\ran}{\ensuremath{\operatorname{ran}}}

\newcommand{\ds}{\ensuremath{\mathrm ds}}

\usepackage{authblk}

\newcommand{\R}{\mathbb{R}}

\DeclareAcronym{nse}{
  short = NSE,
  long  = Navier--Stokes equations
}
\DeclareAcronym{psp}{
  short = SVE,
  long  = Stress-Velocity Estimator 
}

\DeclareAcronym{ps}{
  short = SV,
  long  = stress-velocity  
}

\DeclareAcronym{ppe}{
  short = PPE,
  long  = Pressure Poisson Estimator
}
\DeclareAcronym{ste}{
  short = STE,
  long  = Stokes Estimator
}
\DeclareAcronym{vwerp}{
  short = vWERP,
  long  = Virtual Work-Energy Relative Pressure
}

\DeclareAcronym{fe}{
  short = FE,
  long  = finite element,
  long-plural-form = finite elements
}

\title{A Stress-Based Estimator for Pressure and Stress Recovery from Velocity Measurements}

\author[4,5]{Nicolás Barnafi}
\author[1]{Benjamín Cisternas}
\author[2,3]{Ernesto Castillo}
\author[1,6]{Felipe Galarce\footnote{Corresponding author}}

\affil[1]{School of Civil Engineering. Pontificia Universidad Católica de Valparaíso. Valparaíso, Chile.}
\affil[2]{Universidad de Santiago de Chile, Departamento de Ingeniería Mecánica. Santiago, Chile.}
\affil[3]{Computational Heat and Fluid Flow Lab, Universidad de Santiago de Chile, Santiago, Chile}
\affil[4]{Instituto de Ingeniería Matemática y Computacional \& Facultad de Ciencias Biológicas, Pontificia Universidad Católica de Chile, Chile.}
\affil[5]{Center for Mathematical Modeling, Chile}
\affil[6]{Center for Interdisciplinary Research in Biomedicine, Biotechnology and Well-Being (CID3B). Pontificia Universidad Católica de Valparaíso, Chile.}

\date{}

\begin{document}
\maketitle
\vspace{-2.5cm}

\begin{abstract}
\singlespacing

\noindent\rule{\textwidth}{0.4pt}
\noindent\textbf{Abstract}
\vspace{0.2cm}

\noindent Non-invasive pressure field estimation from velocity measurements is a longstanding engineering problem. We propose, analyze, and test a pressure-recovery method that computes a full stress field from velocity measurements, and leaves the pressure estimation as a cheap post-processing step. The method relies on a stress-velocity first order formulation of the Navier--Stokes equations, and we show that the formulation accounts for deviations from incompressibility in the measured velocity field by construction. In addition, we theoretically establish the convergence of the \ac{fe} approximation scheme, the stability of the stress recovery with respect to finite-resolution velocity measurements, and then validate this theory numerically. Our results show that the proposed estimator is robust in convective flow regimes and remains accurate at reduced spatial resolution, improving upon state-of-the-art pressure-recovery strategies.

\vspace{0.2cm}

\noindent\rule{\textwidth}{0.4pt}
\end{abstract}

\vspace{-0.5cm}
\section{Introduction}


Pressure estimation from velocity measurements is a longstanding challenge in many real-world applications where direct pressure acquisition requires procedures that are too expensive, invasive, or even destructive. This need arises across different areas of fluid mechanics, where velocity fields can often be measured with substantially greater spatial coverage than pressure. In biomedical applications, direct pressure measurements are often invasive and velocity-based approaches have been considered for estimating pressure drops across vascular regions \cite{nolte_2022_review, bertoglio_pdrop, GLM2021} and pressure fields in soft tissues \cite{galarce2023_MRE}.  In pipeline transportation systems, reconstructed pressure fields can be used to identify pressure losses, detect leaks, and monitor transient events such as water hammer phenomena \cite{Brunone2000, galarce2025_bingham}. Similar reconstruction strategies are relevant in turbo-machinery, where particle image velocimetry provides velocity measurements from which pressure distributions around rotating blades can be inferred for performance assessment and design optimization \cite{Liu2006}, and in aerodynamics, where flow-field measurements can be used to estimate surface loading and aerodynamic forces without extensive pressure-sensor instrumentation \cite{vanOudheusden2013}. 

From a mathematical standpoint, the recovery task can be formulated both as an inverse problem \cite{Wilcox2010, tarantola2005} and as a data assimilation procedure \cite{asch2016_dataAss, lahoz2010data, Brunton_Kutz_2019}, where noisy and low-resolution velocity data are enhanced by combining an underlying physical model with the observed data, allowing the estimation of hidden states \cite{GLM2021,GGLM2021}, and ensuring physically coherent predictions. Several methodological approaches have been proposed for pressure recovery, including physics-informed neural-networks \cite{sierpe2025estimationhemodynamicparametersphysics}, manifold learning \cite{Franz16032014}, graph neural networks \cite{romor2026}, and reduced-order variational data assimilation \cite{galarce2023bias,GALARCE2025110374}. 

We focus on a family of methods designed to estimate pressure fields by direct manipulation of the momentum conservation, without offline model pre-training. Among the available alternatives, we find the classical \ac{ppe} \cite{REBHOLZ2020112366}, which is based on taking the divergence of the momentum conservation equation, reducing the field computation to the solution of a Poisson problem with a load vector computed from the velocity data. This method is well-known for underestimating the field peaks when noisy data is considered, mostly due to the artificial smoothness requirements the method imposes on the pressure field \cite{pachecoPdrop}. Two other relevant and recently introduced strategies are the \ac{vwerp} estimator \cite{Marlevi2019}, and the \ac{ste} \cite{svihlova_2016}. Both outperform the \ac{ppe}, and they differ in their robustness to noise, computational cost, and versatility. This family of methods will serve as baselines for the comparisons presented in our numerical examples.

This work proposes a new method based on a stress-velocity formulation of the Navier-Stokes equations \cite{cai2004least, cai2010pseudostress}, whose approximability by means of \ac{fe} is well-established under standard hypotheses \cite{caceres2017mixed,camano2017augmented}. The approach reconstructs the full Cauchy stress tensor from velocity measurements and recovers the pressure as a post-processing step from its isotropic component. The main methodological contribution is therefore the reformulation of pressure recovery as a single-field stress reconstruction problem, rather than as a direct pressure estimation problem. From a practical standpoint, reconstructing the stress also provides direct access to shear-related quantities such as wall shear stress and oscillatory shear index \cite{MELLA2026108145}. From a numerical standpoint, the resulting formulation is elliptic and symmetric, and naturally admits an $\ten H(\dive)$-conforming \ac{fe} discretization, thereby avoiding the auxiliary velocity-pressure saddle-point problems required by alternative approaches such as \ac{ste} and \ac{vwerp}, as well as the artificial boundary conditions commonly required for the correctness of the PPE.

We name our method the \ac{psp}, and we show that it is both well-posed and robust against measurement degradation, including an observability result that couples the stress and the velocity fields, as expected for this kind of formulation. On the computational side, we focus on the resolution of the resulting symmetric variational problem by means of the \ac{fe} method, using Raviart-Thomas (RT) elements, which provide a natural choice for the related infinite-dimensional spaces. We test \ac{psp} in several numerical examples, including convergence tests to verify our theoretical findings, and also assess the robustness of the method in a moderately convective Reynolds-number regime that mimics a blood flow stream.

The paper is organized as follows: In Section \ref{sec:p_recovery} we state the pressure recovery problem along with its mathematical foundations. Section \ref{sec:PSP_model} introduces the Stress-Velocity Estimator and develops its variational formulation, \ac{fe} discretization, and well-posedness analysis. Section \ref{sec:sota} describes two state-of-the-art methods, the \ac{vwerp} and the \ac{ste}, which are used as reference methods for the numerical comparisons. Section \ref{sec:experiments} shows a set of numerical tests to verify the convergence properties of the proposed formulation, assess its sensitivity to measurement resolution and noise, and compare its performance with the reference estimators. Finally, Section 6 summarizes the main conclusions and outlines directions for future work.



\paragraph{Notation:} Throughout this work, we denote  vector and tensor quantities through boldface fonts. We denote with $L^2(\Omega)$ the space of functions $f:\Omega\to\R$ such that $\int|f|^2\,dx<\infty$, with inner product $(f,g) \coloneqq \int fg\,dx$. Vector and tensor functions (with range on $\R^d$ and $\R^{d\times d}$) use boldface $\vec L^2$ and $\ten L^2$. The space $H^1(\Omega)$ (respectively $\vec H^1$, $\ten H^1$) denotes the space of $L^2(\Omega)$ functions $f$ whose weak gradient belongs to $\vec L^2(\Omega)$ as well. Finally, $\vec H(\dive; \Omega)$ denotes the space of vector functions $\vec f:\Omega \to \R^d$ whose divergence $\dive \vec f$ belongs to $L^2(\Omega)$, and $\ten H(\Omega;\dive)$ is the space of $\ten L^2(\Omega)$ tensors with $\vec L^2(\Omega)$ row-wise divergence.

\section{The pressure recovery problem}
\label{sec:p_recovery}

Let us consider a time interval $[0,T]$ and a domain $\Omega \subset \bR^d$, with $d=2$ or $3$. We assume that $\Omega$ is a bounded Lipschitz domain and  decompose the boundary into inlet, outlet, and wall sections as $\partial\Omega= \Gamma =\overline{\Gamma}_{\text{in}}\cup\overline{\Gamma}_{\text{out}}\cup\overline{\Gamma}_{\text{w}}$. We consider the incompressible \ac{nse}:
\begin{equation} \label{eq:NSE}
\left\{
\begin{aligned}
\rho \left( \partial_t \bu
+ \bu \cdot\nabla\bu \right)
- \nabla \cdot (2\mu \nabla^s \bu)
+ \nabla p
 &= \bff_b
&& \text{in } \Omega \times [0,T],\\
\dive \bu &= 0
&& \text{in } \Omega \times [0,T],
\end{aligned}
\right.
\end{equation}
where $\bu:\Omega\times[0,T]\to\bR^d$ is the velocity field, $p:\Omega\times[0,T]\to\bR$ is the pressure field, $\rho \in \bR^+$ is the fluid density, and $\mu$ is the  the dynamic viscosity, either a constant $\mu\in\bR^+$ or a scalar field $\mu:\Omega\times[0,T]\to\bR^+$. The term $\bff_b:\Omega\times[0,T]\to\bR^d$ denotes a body force. The symmetric gradient is defined as $2\nabla^s\bu\coloneqq (\nabla\bu+\nabla\bu^\top)$. A relevant quantity of interest to compute is the pressure drop between the inlet and outlet sections. This quantity is independent of the arbitrary pressure level and is defined as
\begin{equation}
    \delta p \coloneqq \frac{1}{\lvert \Gamma_{\text{in}} \rvert} \int_{\Gamma_{\text{in}}} p \, \ \ds - \frac{1}{\lvert \Gamma_{\text{out}} \rvert} \int_{\Gamma_{\text{out}}} p \,\ \ds. 
    \label{eq:pdrop}
\end{equation}
Some methods directly compute the quantity \eqref{eq:pdrop} from data, such as vWERP, while others such as PPE or STE aim to estimate the full pressure field, from which the pressure drop is subsequently obtained in an inexpensive post-processing step. Our approach also provides the full pressure field, but through a different route: it first reconstructs the stress tensor and then recovers the pressure from its trace in a post-processing step.

Velocity measurements are idealized here and generated synthetically from ground-truth velocity fields using two different representations. The simpler approach consists of interpolating the velocity field onto the \ac{fe} meshes with spatial resolutions comparable to that of the  measurements \cite{bertoglio_pdrop}. Alternatively, the measurements can be understood as a set of linear functionals acting on the ground truth field (typically with $\cH = \vec H^1(\Omega)$), as in \cite{HAIK2023115868, cohen2022_nonlinearSpaces, GLM2022}:
\begin{equation}
\ell_i(\bv)
=
\left( \bv, \boldsymbol{\omega}_i \right)_\cH,
\qquad
\forall \bv \in \cH,
\end{equation}
for $i=1,\ldots,m$, where $\boldsymbol{\omega}_i\in\cH$ is the Riesz representer  associated with
the $i$-th sensor $\ell_i\in\cH'$, where $\cH'$ denotes the dual space of $\cH$. We can therefore interpret the degraded velocity
field  $\bu_m$ as an element of the observation space
$W_m=\operatorname{span}\{\boldsymbol{\omega}_1,\ldots,\boldsymbol{\omega}_m\}
\subset\cH$, and we can generate synthetic velocity measurements by applying the projection operator
onto $W_m$ to the ground truth fields: $ \bu_m = \Pi_m\bu \in W_m$, where $\Pi_m:\cH\to W_m$ denotes the $\cH$-orthogonal projector onto $W_m$. Thus, $\bu_m$
represents the velocity information available for the pressure-recovery problem.
Depending on the support and spatial distribution of the representers
$\boldsymbol{\omega}_i$, the operator $\Pi_m$ can model spatial under-resolution,
voxel-wise measurements, or observations restricted to a subset of the domain. In
the numerical tests presented in Section~\ref{sec:experiments}, the observation space $W_m$ is induced by voxel grids of
different resolutions.

Although the ideal velocity field $\bu$ is divergence-free, the measured field
$\bu_m$ may not satisfy this property exactly due to noise, partial measurements, or even modeling errors.
We show that the proposed formulation associates the recovered pressure with a corrected divergence-free velocity field through a bias term that accounts for deviations in both the incompressibility and momentum balance equations in the measured data. In other methods, a solenoidal projection could be introduced as a preprocessing step. However, when
pressure-recovery estimators are driven by measured velocity data \cite{svihlova_2016, bertoglio_pdrop} rather than by an exact \ac{nse} velocity field,
their outputs should be understood as pressure estimates associated with the available measurement field, rather than as exact pressure reconstructions corresponding to a fully consistent incompressible flow.

\noindent The variational form \eqref{eq:ps_variational_form}, together with the boundary conditions define the continuous \acl{ps} problem. In this equivalent reformulation of the \ac{nse}, the stress tensor $\ten\sigma$ is the field to be reconstructed, while all terms depending on the velocity are evaluated from the known field $\bu$. After $\ten\sigma$ is computed, the pressure is recovered independently from \eqref{eq:pressure_trace}. Thus, the pressure recovery problem is reduced to a stress-reconstruction problem driven entirely by the velocity data. In the following section, this continuous formulation is discretized in time and space, leading to the \ac{fe} problem used in the numerical implementation.

\section{\acl{psp} }
\label{sec:PSP_model}
In this section, we introduce the proposed \ac{psp}. The method is based on reconstructing a stress field from a known velocity field, which can be obtained from a numerical simulation or from experimental measurements. The objective is to recover pressure from the reconstructed stress tensor, rather than computing the pressure directly by solving a Poisson-type equation or an auxiliary Stokes problem.

The key idea is to replace the original pressure-recovery problem with the reconstruction of a stress field that satisfies two complementary equations: the momentum conservation equation and a constitutive equation. The \acl{ps} structure is introduced through the constitutive relation, which links the deviatoric component of the stress tensor to the viscous stress. Once the stress field is obtained, the pressure is recovered from its trace in a post-processing step. In this section we present the continuous formulation, its \ac{fe} discretization, and the well-posedness analysis of the proposed formulation.

\subsection{Continuous formulation}

At the continuous level, starting from the momentum equation in the \ac{nse} \eqref{eq:NSE}, we rewrite the balance by keeping the surface-force contribution on the left-hand side and grouping the inertial and body-force terms on the right-hand side. This yields:
\begin{equation}
    -\nabla p + \nabla \cdot (2\mu \nabla^s \bu) = \rho \left( \partial_t \bu + \bu \cdot \nabla \bu \right)
- \bff_b.
\label{eq:NSE_momentum}
\end{equation}

We define the right-hand side as
$\bff(\bu) \coloneqq \rho\left(\partial_t\bu +\bu\cdot\nabla\bu\right)-\bff_b$, which contains the inertial and body-force contributions. We also define the Cauchy stress tensor $\ts:\Omega\times[0,T]\to\bR^{d\times d}$, which accounts for the pressure and viscous contributions to the surface traction. It is defined by the constitutive equation:
\begin{equation}
    \ts \coloneqq  -p \ten{I}_d+\ten \tau_v(\bu) \qquad \text{in } \Omega\times[0,T],
\label{eq:constitutive}
\end{equation}
where $\ten \tau_v(\bu)=2\mu \nabla^s \bu$ is the viscous stress tensor. No additional assumptions regarding the explicit rheological law are required for the pressure recovery. With these definitions, the momentum equation can be written in terms of the stress tensor as:
\begin{equation}
\label{eq:ps_momentum}
\nabla\cdot\ten\sigma=\bff(\bu)\qquad
\text{in } \Omega\times[0,T].
\end{equation}
Notice that, in this system, the pressure remains an explicit variable. A direct pressure-recovery formulation would require eliminating the stress tensor, leading to a Poisson-type pressure equation. Instead, we treat the stress tensor as the field to be reconstructed and recover the pressure a posteriori from its trace. For this purpose, we introduce the deviatoric part of the stress, defined as
$\ts^d\coloneqq \ts-\frac{1}{d}\tr(\ts)\ten I_d$. Since the pressure contribution corresponds to the isotropic part of the Cauchy stress tensor, the deviatoric projection eliminates the explicit pressure term from \eqref{eq:constitutive}, which becomes
\begin{equation}
    \ts^d = \ten \tau_v(\bu).
    \label{eq:ps_deviatoric_relation}
\end{equation}
This deviatoric relation does not reduce the recovered field to the deviatoric stress alone. The formulation still recovers the full Cauchy stress tensor, since the momentum equation couples the complete tensor through $\dive\ts$, while the deviatoric relation imposes the constitutive constraint on its non-isotropic part. Thus, the constitutive model is enforced in deviatoric form without discarding the isotropic contribution. Finally, after computing the stress tensor, the pressure field is recovered from the trace of \eqref{eq:constitutive}. Since $\ten\tau_v(\bu)$ is traceless, this gives
\begin{equation}
p = -\frac{1}{d}\tr(\ts).
\label{eq:pressure_trace}
\end{equation}

We now derive the variational form associated with
\eqref{eq:ps_momentum} and \eqref{eq:ps_deviatoric_relation}. We first consider the
momentum equation separately. Since this equation is written in terms of the divergence of the stress tensor and the velocity field is known, the test function at this stage is required only to belong to $\vec L^2(\Omega)$. Therefore, the weak form of the momentum equation is given by
\begin{equation}
(\dive\ten\sigma,\vec v)
=
( \bff(\bu),\vec v)
\qquad
\forall \vec v\in\vec L^2(\Omega).
\label{eq:weak_momentum}
\end{equation}
The weak form of the constitutive relation is obtained by testing \eqref{eq:ps_deviatoric_relation} with the deviatoric part of the tensor test function. Thus,
\begin{equation}
(\ten\sigma^d,\ten\tau^d)
=
(\ten \tau_v(\bu),\ten\tau^d)
\qquad
\forall \ten\tau\in\ten H(\dive;\Omega).
\label{eq:weak_deviatoric}
\end{equation}
Since the deviatoric operator is a linear and bounded projection in $\ten L^2(\Omega)$, for every $\ten\tau\in\ten H(\dive;\Omega)$ the tensor $\ten\tau^d$ is admissible as a test function in \eqref{eq:weak_deviatoric}. Having defined the weak forms of both relations, we choose $\vec v=\dive\ten\tau$ in \eqref{eq:weak_momentum}. Combining the resulting equation with
\eqref{eq:weak_deviatoric} gives the following variational form for the \acl{ps} recovery problem: find $\ten\sigma\in\ten H(\dive;\Omega)$ such that

\begin{equation}
\left\{
\begin{aligned}
(\dive\ten\sigma,\dive\ten\tau)
+
(\ten\sigma^d,\ten\tau^d)
&=
(\ten\tau_v(\bu),\ten\tau^d)
+
(\bff(\bu),\dive\ten\tau)
&&
\forall \ten\tau\in\ten H(\dive;\Omega),
\\
\ten\sigma\vec n
&=
\ten\sigma_D\vec n,
&&
\text{on } \Gamma_{\text{D}} .
\end{aligned}
\right.
\label{eq:ps_variational_form}
\end{equation}
The boundary conditions in the stress-velocity formulation differs from that of the standard velocity-pressure formulation because it is a mixed formulation. This implies that the roles of Dirichlet and Neumann boundary conditions are flipped, and thus Neumann boundary conditions become essential. In particular, this means that no-slip boundary conditions are represented simply by omitting boundary terms in the corresponding boundaries, and the essential conditions prescribed abstractly above correspond to a normal traction given by $\ten \sigma_D = \left( -p_D \ten I + \ten\tau_v(\vec u)\right) \vec n$.

\subsection{Correction for non-solenoidal measurements}

One relevant property of the proposed formulation is that the recovered stress tensor is associated with a divergence-free velocity field, even when the measured velocity field is not. For simplicity, we will carry out our analysis by assuming that $\ten \tau_v(\bu) = \grad \bu$.  To examine this property, we consider the full-measurement case $\bu_m = \bu$ and introduce the space $Z \coloneqq \{ \ten \tau \in  \ten H_0(\dive; \Omega) : \dive \ten \tau = 0 \}$, with $\ten H_0(\dive; \Omega) \coloneqq \{ \ten \tau \in \ten H(\dive; \Omega) : \ten \tau \cdot \vec n = 0 \text{ on } \Gamma_{\text{D}} \}$. Let us take $\ten \tau \in Z$ in \eqref{eq:ps_variational_form}, which leads to: 
\begin{equation}
(\ten \sigma_d - \ten \tau_v(\vec u), \ten \tau ) = 0 \quad \forall \ten \tau \in Z,
\label{eq:weak_projZ}
\end{equation}
where $(\ten \sigma^d, \ten \tau) = (\ten \sigma^d, \ten \tau^d)$, because $\ten \sigma^d$ is traceless. We note that this implies that $\ten G\coloneqq \ten \sigma_d - \ten \tau_v(\vec u)$ is orthogonal to $\ker \dive$, and thus the Closed Range Theorem yields that $\ten G$ belongs to $\left( \ker \dive \right)^\perp = \ran\left(\grad\right)$, so that there exists some vector function $\bz$ in $\vec H^1(\Omega)$ with homogeneous boundary condition on $\Gamma_N$  such that $\ten G =  \nabla \bz$ \cite{girault2012finite}. This yields a biased constitutive law: 
\begin{equation}
\ten \sigma^d = \grad \bu +  \nabla \bz.
\label{eq:biased_cons}
\end{equation}
Taking the trace in \eqref{eq:biased_cons} then gives:
\begin{equation}
    \dive \bu = - \dive \bz,
\end{equation}
so that the recovered stress always corresponds with an incompressible velocity field $\vec w = \vec u + \vec z$. Furthermore, if the measured velocity field is divergence-free, no correction is required and $\bz=0$.

The correction also enters the momentum balance. To see this, consider the explicit pseudo-stress Newtonian formulation for the corrected deviatoric stress $\ten \sigma^d = \ten \tau_v = \nabla \bw$, meaning that \eqref{eq:ps_variational_form} gives:
\begin{align*}
(\dive\ten\sigma, \dive\ten\tau)
+ ( \grad (\bu + \bz),\ten\tau)
& =
(\grad \bu, \ten\tau) +
(\bff(\bu),\dive\ten\tau)
\quad \forall \ten\tau\in\ten H(\dive;\Omega),
\end{align*}
which yields after integration by parts:
\begin{align*}
(\dive\ten\sigma - \bz - \bff(\bu),\dive \ten\tau)
& = 0 \quad \forall \ten\tau\in\ten H(\dive;\Omega).
\end{align*}
It remains to establish that the strong residual vanishes, rather than being merely orthogonal to the range of the divergence operator. Since $\dive:\ten H(\dive; \Omega) \to \vec L^2(\Omega)$ is surjective \cite{BrezziFortin1991} (to $L_0^2(\Omega)$ in the pure Dirichlet case), then testing against all possible variations in $\vec L^2(\Omega)$ is equivalent to $\dive \ten \tau$ for $\ten \tau$ in $\ten H (\dive; \Omega)$. Thus:
$$
\dive \ten \sigma = \bff(\bu) + \bz.
$$
The resulting corrected first-order system can therefore be written as
\begin{align*}
    \dive \ten \sigma &= \bff(\bu) + \bz, \\
    \ten \sigma^d &= \nabla (\bu + \bz). \\
\end{align*}

Thus, this modified formulation recovers the stress field associated with the corrected force--velocity pair $(\bff(\bu) + \vec z, \bu + \vec z)$.  Furthermore, when the measured velocity fields satisfy $\dive \vec u = 0$, no correction is introduced, and $\vec z =0$, i.e., the data is compatible with the strong formulation. Conversely, if the data are not compatible, the discrepancy is represented by the correction field $\vec z$ that one can compute numerically. We highlight that if we consider a different rheology, the computations remain the same, but the interpretation of $\bz$ is simply less clear as a velocity correction, so we have preferred this simpler scenario for the analysis.

\subsection{Discrete formulation}

\noindent To obtain the \ac{fe} approximation of
\eqref{eq:ps_variational_form}, let $\mathcal T_h=\cup_i K_i$ be a triangulation of $\Omega$. For an integer $k\geq 0$, we consider the Raviart--Thomas \ac{fe} space of degree $k$, denoted by $\mathcal{RT}^k$, and define the $\ten H(\dive;\Omega)$-conforming discrete space
$$ \ten H_h = \{\ten \tau\in\ten L^2(\Omega): \ten\tau_{(i, \cdot)}|_K \in \mathcal{RT}^k, K\in \mathcal T_h, i\in\{1,\dots,d\}\}. $$
This tensor-valued space is constructed row-wise, meaning that each row of the stress tensor is approximated as a vector field in a Raviart--Thomas space. The discrete problem is: find $\ten\sigma_h\in\ten H_h$ such that
\begin{equation}
(\dive\ten\sigma_h,\dive\ten\tau_h)
+
(\ten\sigma_h^d,\ten\tau_h^d)
=
(\ten \tau_v(\bu),\ten\tau_h^d)
+
(\vec f(\bu),\dive\ten\tau_h)
\qquad
\forall \ten\tau_h\in\ten H_h.
\label{eq:ps-discrete}
\end{equation}

\noindent The formulation in \eqref{eq:ps-discrete} provides the \ac{fe} spatial approximation of the continuous \acl{ps} problem at a fixed time. In transient flows, however, the right-hand side also depends on the material acceleration of the known velocity field. Therefore, the time derivative in $\vec f(\bu)$ must be discretized to evaluate it from the available sequence of velocity fields.

Let $t_n=n\Delta t$, with $\Delta t>0$ denoting the time step. We assume that the measured velocity field is available at the discrete times $t_n$. For $n\geq 1$, the time derivative is approximated using the second-order backward differentiation formula (BDF2),
\begin{equation*}
\partial_t \bu(t_{n+1})
\approx
\frac{
3\bu^{n+1}
-
4\bu^{n}
+
\bu^{n-1}
}{2\Delta t}.
\label{eq:BDF2_velocity}
\end{equation*}

\noindent
Therefore, the right-hand side of the momentum equation at time $t_{n+1}$ is evaluated as
\begin{equation}
\vec f^{\,n+1}(\bu)
=
\rho\left(
\frac{
3\bu^{n+1}
-
4\bu^{n}
+
\bu^{n-1}
}{2\Delta t}
+
\bu^{n+1} \cdot \nabla\bu^{n+1}\right) -\bff_b^{\,n+1}.
\label{eq:discrete_rhs_BDF2}
\end{equation}
At each time step $t_{n+1}$, the fully discrete \acl{ps} problem reads: find $\ten\sigma_h^{n+1}\in \ten H_h$ such that
\begin{equation}
(\dive\ten\sigma_h^{n+1},\dive\ten\tau_h)
+
((\ten\sigma_h^{n+1})^d,\ten\tau_h^d)
=
(\ten \tau_v(\bu^{n+1}),\ten\tau_h^d)
+
( \vec f^{\,n+1}(\bu),\dive\ten\tau_h)  \qquad \forall \ten\tau_h\in\ten H_h.
\label{eq:ps-discrete-BDF2}
\end{equation}
\noindent The pressure field at $t_{n+1}$ is then obtained from the post-processing  equation (\ref{eq:pressure_trace}). Hence, the fully discrete \ac{psp} reduces, at each time step, to a stress-reconstruction problem driven by the known velocity field.

\subsection{Well‑posedness analysis}
In this section, we establish the well-posedness of problem \eqref{eq:ps_variational_form} by a standard application  of the Lax-Milgram Lemma \cite{EG2013}. Since the discrete space $\ten H_h$ is conforming, the same framework can subsequently be used to analyze the discrete problem. To this end, consider the bilinear form associated with\eqref{eq:ps_variational_form}, given  by $a:\ten H(\dive; \Omega)\times \ten H(\dive; \Omega) \to \R$ such that
        $$ a(\ten \sigma, \ten \tau) = (\dive \ten \sigma, \dive \ten \tau) + (\ten \sigma^d, \ten \tau^d). $$
Boundedness follows from the continuity of the deviatoric operator, $\| \ten \sigma^d \|_{L^2(\Omega)}\leq C  \| \ten \sigma \|_{L^2(\Omega)}$:
        $$ a(\ten \sigma, \ten \tau) \leq \|\dive \ten \sigma\|_{L^2(\Omega)} \|\dive \ten \tau\|_{L^2(\Omega)} + C^2 \| \ten \sigma \|_{L^2(\Omega)} \| \ten \tau \|_{L^2(\Omega)} \leq \max\{1, C^2\} \|\ten \sigma\|_{H(\dive;\Omega)} \| \ten \tau \|_{H(\dive;\Omega)} $$
For coercivity, we use the estimate established in \cite{arnold1984family}, according to which, under the appropriate boundary or normalization condition, there exists  $c_1>0$ such that
    $$ c_1\|\ten \tau \|_0^2 \leq \| \ten \tau^d \|_0^2 + \| \dive \ten \tau \|_0^2. $$
This result yields the ellipticity of $a$ as follows:
\begin{multline*}
a(\ten \tau, \ten \tau) = \| \ten \tau^d \|^2_{L^2(\Omega)} + \| \dive{\ten \tau} \|^2_{L^2(\Omega)} \\ 
    = \frac 1 2 \left(\|\ten \tau^d\|_0^2+\|\dive \ten \tau\|_0^2\right) + \frac 1 2 \left(\|\ten \tau^d\|_0^2 + \|\dive \ten \tau \|_0^2\right) \geq \min\{1/2,c_1\} \| \ten \tau \|_{H(\dive;\Omega)}^2.
\end{multline*}
The continuity and coercivity of $a$, together with the continuity of the linear functional on the right-hand side of \eqref{eq:ps_variational_form}, satisfy the hypotheses of the Lax--Milgram lemma, which guarantees the existence and uniqueness of a solution to the stress-reconstruction problem. 

\paragraph{Convergence of the discrete scheme:} Since the continuous problem is well-posed by Lax-Milgram lemma, standard \ac{fe} theory ensures that any conforming approximation space $\ten H_h\subseteq \ten H$ yields a well-posed and convergent discrete problem. For the chosen Raviart--Thomas discretizations, the convergence rate is $O(h)$ \cite{EG2013}.

\paragraph{Convergence of measurements:} To analyze the convergence of the method with respect to finite measurements, we consider the reconstruction obtained from partial observations, $\ten \sigma_m$, and the corresponding exact reconstruction, $\ten \sigma$. Subtracting their governing equations yields the following equation for the error $\ten e_{\sigma,m}\coloneqq \sigma - \sigma_m$:
$$ ( \dive{\ten e_{\sigma,m}}, \dive{\ten \tau}) + (\ten e_{\sigma,m}^d, \ten \tau^d) = (\tau_v (\vec u) - \tau_v(\vec u_m), \ten \tau^d) +  (\vec f(\vec u) - \vec f(\vec u_m), \dive \ten \tau). $$
We consider the measurement error $\vec e_{u,m} \coloneqq \vec u - \vec u_m$. The contribution of each term is bounded as follows:
$$ \|\tau_v(\vec u) - \tau_v(\vec u_m)\| \leq C \| \vec e_{u,m} \|_{\vec H^1(\Omega)}, $$
using $C$ as a placeholder for a generic positive constant $C$:
$$ \begin{aligned} 
    \| \vec f(\vec u) - \vec f(\vec u_m) \| &= \| \partial_t \vec e_{u,m} + \vec u\cdot \nabla \vec u - \vec u_m \cdot \nabla \vec u_m \| \\ 
                        &\leq \| \partial_t \vec e_{u,m} \|_{L^2(\Omega)} + C \|\vec e_{u,m} \|_{H^1(\Omega)} \\
                        &\leq M \left( \| \partial_t \vec e_{u,m} \|_{L^2(\Omega)} + \| \vec e_{u,m} \|_{H^1(\Omega)} \right).
    \end{aligned} $$
    We have denoted by $M$ some positive constant depending on the norms of the solution $\vec u$. This shows that the recovery error can be estimated by
    \begin{equation}
         \| \ten e_{\sigma, m} \|_{H(\dive;\Omega)} \leq M \left( \| \partial_t \vec e_{u,m} \|_{L^2(\Omega)} + \| \vec e_{u,m} \|_{H^1(\Omega)} \right), 
    \label{eq:recovery_error}
    \end{equation}
        which implies that if the measurement operator converges in $L^2(0,T;H^1(\Omega))$, then convergence is guaranteed. 
\noindent As a wrap-up of this section, we have established that: 
    \begin{itemize}
        \item The \ac{psp} formulation corresponds to the recovery of an incompressible flow, and the data mismatch can be computed explicitly.
        \item The \ac{psp} formulation is well-posed and the proposed discrete scheme is convergent.
        \item The infinite measurements scenario can be approximated through partial measurements, and convergence is therein guaranteed. 
    \end{itemize}

\section{Baseline methods}
\label{sec:sota}

This section reviews other state-of-the-art methods used to compare the proposed \ac{psp} against. We focus on the \acl{vwerp} \cite{bertoglio_pdrop, Marlevi2019, marlevi21} and the \acl{ste} \cite{svihlova_2016}, since both methods have reported better accuracy and robustness against noisy and noise-free measurements than the classical \acl{ppe}. For completeness, we also recall the PPE formulation as the standard baseline approach for pressure recovery from velocity data. We briefly describe these approaches, highlighting their main assumptions and implementation details.

\subsection{Pressure Poisson Estimator}
The \ac{ppe} is based on solving a Poisson equation for the pressure, whose right-hand side is computed from the inertial terms of the \ac{nse}. Starting from the momentum equation and taking its divergence, we obtain
\begin{equation}
-\Delta p
=
\nabla \cdot
\left(
\rho \partial_t \vec u
+
\rho(\vec u\cdot\nabla)\vec u
-
\mu \Delta \vec u
\right).
\label{eq:PPE}
\end{equation}
For an incompressible Newtonian fluid, the viscous contribution vanishes at the continuous level, since $\nabla\cdot\Delta\vec u=\Delta(\nabla\cdot\vec u)=0$. Therefore, the Laplacian of the pressure is driven by the divergence of the inertial terms. This leads to the following weak form:
\begin{equation} \label{eq:weak_PPE}
    (\nabla p, \nabla v) = -\rho(\partial_t \vec u,\nabla v) - \rho (\vec u \cdot \nabla \vec u,\nabla v)  \qquad \forall v \in H^1 (\Omega).
\end{equation}
In practice, for the spatial discretization,  piecewise linear elements are used, while the time derivative can be approximated using a first or second order backward finite difference. The pressure level is fixed by imposing suitable boundary conditions, following the standard implementation of the method. 

This estimator recovers the pressure field from a known velocity field. Once the pressure field has been computed, the pressure drop can be evaluated using \eqref{eq:pdrop}. Although the \ac{ppe} formulation is the standard baseline approach for pressure recovery from velocity data, Bertoglio et al. \cite{bertoglio_pdrop} show that alternative estimators, such as \ac{vwerp} and \ac{ste}, can provide improved accuracy and robustness, especially in the presence of noisy velocity measurements.

\subsection{Virtual Work-Energy Relative Pressure Estimator}

The \ac{vwerp} estimator \cite{bertoglio_pdrop, Marlevi2019, marlevi21} is based on the weak form of the \ac{nse} for an incompressible Newtonian fluid. In this approach, an auxiliary test function is computed in advance and chosen so that the weak pressure term depends only on the inlet and outlet sections. As a result, the pressure drop can be estimated without reconstructing the complete pressure field. More precisely, the test function $\vec v \in \vec H^1(\Omega)$ is chosen to satisfy 
$\nabla \cdot \vec v = 0  \text{ in } \Omega,$ and $ \vec v \cdot \vec n = 0  \text{ on } \Gamma_w.$

\noindent In practice, this field is obtained from a Stokes auxiliary variational problem, and remains fixed during the pressure-drop evaluation. The weak form of the momentum equation in \eqref{eq:NSE} is then derived using $\vec v$ as test function. Under the incompressibility constraint and assuming constant viscosity, the viscous term
reduces to $\nabla\cdot(2\mu\nabla^s\bu)=\mu\Delta\bu$. The convective, pressure, and viscous terms are integrated by parts to isolate the boundary contributions associated with the pressure drop:

\begin{equation}
    I_{kin}(\bu) = \rho \int_{\Omega}  \partial_t \bu \cdot \vec{v} \ dx,
\end{equation}
\begin{equation}
    I_{conv}(\textbf{u}) = -\rho \int_{\Omega} (\textbf{u} \cdot \nabla \textbf{v}) \cdot \textbf{u} \ dx + \rho \int_{\Gamma} (\textbf{u} \cdot \textbf{n})(\textbf{u} \cdot \textbf{v}) \ ds,
\end{equation}
\begin{equation}
    I_{pres}(\textbf{u}) = -\int_{\Omega} p (\nabla \cdot \textbf{v}) \ dx +\int_{\Gamma} p(\textbf{v} \cdot \textbf{n}) \ ds = \int_{\Gamma_i \cup \Gamma_o} p(\textbf{v} \cdot \textbf{n}) \ ds,
\end{equation}
\begin{equation}
    I_{vis}(\textbf{u}) =  \mu\int_{\Omega} \nabla \textbf{u} : \nabla \textbf{v} \ dx -\mu\int_{\Gamma} (\nabla \vec{u} \cdot \textbf{n}) \cdot \textbf{v} \ ds.
\end{equation}

\noindent Since the kinetic, convective, and viscous terms depend only on the known velocity field and on the precomputed test function $\vec v$, they can be evaluated directly. Moreover, assuming that the pressure is approximately constant on the inlet and outlet sections, the pressure drop can be obtained from the pressure term as
\begin{equation}
    I_{\mathrm{pres}}(\bu) = p_i \int_{\Gamma_i} \vec v \cdot \vec n   \ ds + p_o \int_{\Gamma_o} \vec v \cdot \vec n  \ ds .
\end{equation}
Since $\vec v$ is divergence-free in $\Omega$ and satisfies $\vec v\cdot\vec n=0$ on $\Gamma_w$, the boundary flux is balanced between the inlet and outlet sections. Therefore, the flux through the inlet is equal to the flux through the outlet with opposite sign. Defining
\begin{equation}
\delta p = p_i-p_o, \qquad \Lambda(\vec v)
= \int_{\Gamma_i} \vec v\cdot\vec n \, \ ds.
\end{equation}
The pressure term can be written as $I_{\mathrm{pres}}(\bu) = \delta p\,\Lambda(\vec v)$. Finally the pressure drop is given by:

\begin{equation}
    \delta p (\vec u) = -\frac{1}{\Lambda(\vec{v})}\left(I_{kin}(\partial_t \vec u)+I_{conv}(\vec u)+I_{vis}(\vec u)\right).
    \label{eq:vwerp}
\end{equation}
In this way, equation \eqref{eq:vwerp} computes the pressure drop directly from velocity-dependent integral terms and a precomputed auxiliary test function. This makes the method computationally efficient, but it does not reconstruct the pressure field. The output is limited to a global pressure-drop estimate under the assumption of nearly constant pressure over the inlet and outlet sections.

\subsection{Stokes Estimator}
The \ac{ste} \cite{svihlova_2016} is based on solving an auxiliary Stokes problem whose right-hand side is computed from the inertial and viscous terms of the \ac{nse}. Given the known velocity field, the method computes an auxiliary field $\vec w$ and a pressure field such that the momentum residual associated with $\bu$ is balanced by a Stokes-type system. This leads to the following auxiliary Stokes problem:
\begin{equation} \label{eq:STE}
\left\{
\begin{aligned}
-\Delta \textbf{w} + \nabla p = \vec f(\vec{u})
&& \text{in } \Omega, \\
\nabla \cdot \textbf{w} = 0
&& \text{in } \Omega, \\
\textbf{w} = 0
 && \text{on } \Gamma,
\end{aligned}
\right.
\end{equation}
where the right-hand side for an incompressible Newtonian fluid is given by $ \vec f(\vec u) = -(\rho \partial_t \vec u + \rho \vec u \cdot \nabla \vec u) + \mu \Delta \vec u$. Notice that the auxiliary field is required to satisfy $\vec w=\vec 0$ on the whole boundary. Therefore, a suitable pressure condition must be imposed to fix the arbitrary pressure level and obtain a uniquely defined pressure solution.

The robustness of the STE formulation can be improved by integrating by parts \cite{bertoglio_pdrop}. This reduces the direct differentiation of noisy data and leads to a more stable evaluation of the right-hand side. The corresponding momentum equation reads
\begin{equation}
(\nabla \vec w, \nabla \vec v)
-(p, \nabla\cdot \vec v)=
-{\rho}(\partial_t \vec u, \vec v)
+\rho(\vec u\cdot\nabla\vec v,\vec u)
-\mu(\nabla \vec u, \nabla \vec v) \quad \forall \vec v \in \vec H^1 (\Omega).
\label{eq:STEint}
\end{equation}
In practice, Taylor--Hood elements are used for the spatial discretization,  while the time derivative can be approximated using a first or second order backward finite difference. As in the PPE, this estimator recovers the pressure field from a known velocity field. Once the pressure field has been computed, the pressure drop can be evaluated using \eqref{eq:pdrop}.

\section{Numerical tests}
\label{sec:experiments}
In this section, we present a set of two-dimensional flow numerical tests designed to validate the proposed methodology. The numerical tests were implemented using the open-source \ac{fe} library Firedrake \cite{FiredrakeUserManual}, and the resulting linear
systems were solved with a direct LU factorization using the Multifrontal Massively Parallel Solver (MUMPS) \cite{MUMPS:1,MUMPS:2}. The experiments are organized to assess different aspects of the method:

\begin{enumerate}
    \item A convergence test for the \acl{ps} recovery operator using an analytical stress tensor, used to verify the convergence of the \ac{fe} scheme.

    \item A velocity-measurement quality test, used to evaluate the influence of the spatial quality of the input velocity measurements on the reconstruction accuracy and convergence rate of the recovered stress field.

    \item A Venturi channel benchmark, used to compare the proposed \ac{psp}  with state-of-the-art pressure recovery estimators, including vWERP \cite{bertoglio_pdrop} and STE \cite{svihlova_2016}.

    \item A noisy-data Venturi channel benchmark, used to assess the robustness of the proposed \ac{psp} estimator under perturbed velocity measurements and to compare its performance with the same reference estimators.
\end{enumerate}

\subsection{Convergence test}
To assess the spatial convergence of the \ac{fe} approximation, we first consider the \acl{ps} recovery operator independently of the pressure reconstruction problem, focusing only on the bilinear form associated with the \ac{psp} formulation. In this test, the right-hand side is prescribed directly, so that no velocity field is required.

The test uses a manufactured stress tensor to compute the \ac{fe} approximation error. The right-hand side is chosen so that a prescribed analytical stress field is the exact solution of the corresponding strong problem \eqref{eq:strong_pseudostress_recovery}, which is obtained from the elliptic differential operator associated with the bilinear form \eqref{eq:ps_variational_form}.
\begin{equation}
\label{eq:strong_pseudostress_recovery}
\begin{cases}
\ten{\sigma}^{\mathrm d}
-
\nabla(\dive \ten{\sigma})
=
\ten{f},
& \text{in } \Omega, \\
\ten{\sigma} \vec n
=
\ten{\sigma}_{\mathrm D} \vec n,
& \text{on } \Gamma .
\end{cases}
\end{equation}
The analytical stress field used as the reference solution is prescribed as
\begin{equation}
\label{eq:sigma_exact}
\ten{\sigma}(x,y)
=
\begin{pmatrix}
\sin(x)\cos(y) & x^3+y^3 \\[2pt]
x^3+y^3 & \sin(x)\cos(y)
\end{pmatrix}.
\end{equation}
This stress field contains both isotropic and deviatoric contributions, and therefore provides a nontrivial manufactured solution for the \acl{ps} recovery problem. The corresponding source term $\ten f$ is obtained by applying the strong operator \eqref{eq:strong_pseudostress_recovery} to \eqref{eq:sigma_exact}. Its components are given by
\begin{equation}
\label{eq:rhs_sigma_components}
\begin{cases}
f_{xx} = f_{yy} = \sin(x)\cos(y),
\\
f_{xy}= x^3+y^3+\cos(x)\sin(y)-6y,
\\
f_{yx}=x^3+y^3+\cos(x)\sin(y)-6x.
\end{cases}
\end{equation}
The convergence test is performed on a rectangular domain discretized with triangular meshes obtained from successive uniform refinements. A total of eight meshes are considered, ranging from 64 elements in the coarsest configuration to $658{,}432$ elements in the finest one. We denote by $h$ the characteristic mesh size, taken as the maximum diameter of the elements in the mesh. The stress tensor is approximated using lowest-order Raviart--Thomas elements. For each mesh, the numerical solution $\ten\sigma_h$ is compared with the analytical
stress tensor $\ten\sigma$. We report the errors
\begin{equation}
e_{\ten L^2}(h)
=
\|\ten\sigma-\ten\sigma_h\|_{\ten L^2(\Omega)},
\qquad
e_{\ten H(\dive)}(h)
=
\|\ten\sigma-\ten\sigma_h\|_{\ten H(\dive;\Omega)}.
\end{equation}
The corresponding convergence rates are obtained from successive mesh refinements.
The results are summarized in Table~\ref{tab:sigma_spatial_convergence}. The observed first-order decay of the errors is consistent with the lowest-order Raviart--Thomas approximation, supporting the consistency of the proposed \acl{ps} recovery operator and its \ac{fe} discretization.

\begin{table}[htbp]
\centering
\caption{Spatial convergence history for the \acl{ps} recovery operator.}
\begin{tabular}{ccccc}
\toprule
$h$ & $\|\ten \sigma-\ten \sigma_h\|_{\ten L^2(\Omega)}$ & rate
& $\|\ten \sigma- \ten \sigma_h\|_{\ten H(\mathrm{div};\Omega)}$ & rate \\
\midrule
$1.414\times10^{-2}$ & $1.832\times10^{-4}$ & -
& $1.876\times10^{-4}$ & - \\
$7.071\times10^{-3}$ & $9.160\times10^{-5}$ & $1.000$
& $9.378\times10^{-5}$ & $1.000$ \\
$3.536\times10^{-3}$ & $4.580\times10^{-5}$ & $1.000$
& $4.689\times10^{-5}$ & $1.000$ \\
$1.768\times10^{-3}$ & $2.290\times10^{-5}$ & $1.000$
& $2.345\times10^{-5}$ & $1.000$ \\
$8.839\times10^{-4}$ & $1.145\times10^{-5}$ & $1.000$
& $1.172\times10^{-5}$ & $1.000$ \\
$4.419\times10^{-4}$ & $5.725\times10^{-6}$ & $1.000$
& $5.862\times10^{-6}$ & $1.000$ \\
$2.210\times10^{-4}$ & $2.862\times10^{-6}$ & $1.000$
& $2.931\times10^{-6}$ & $1.000$ \\
$1.105\times10^{-4}$ & $1.431\times10^{-6}$ & $1.000$
& $1.465\times10^{-6}$ & $1.000$ \\
\bottomrule
\end{tabular}
\label{tab:sigma_spatial_convergence}
\end{table}

\subsection{Velocity-measurement convergence test}
We now consider a steady manufactured flow to assess the convergence of the recovered stress field with respect to the spatial resolution of the velocity measurements. The analytical velocity field is sampled over a sequence of voxel grids, producing voxelized velocity fields that are used as input data in the \ac{psp}. The reference solution used for the error evaluation satisfies a Stokes-type problem with a prescribed body force on the right-hand side. The analytical velocity and pressure fields are given by

\begin{equation}
\begin{aligned}
\vec u(x,y)
&=
\begin{pmatrix}
\pi\cos(\pi y)\left(\sin(\pi x)+\dfrac{1}{4}\sin(2\pi x)\right)
\\
-\pi\sin(\pi y)\left(\cos(\pi x)+\dfrac{1}{2}\cos(2\pi x)\right)
\end{pmatrix},
\qquad
p(x,y)=0.55-0.9x-0.2y.
\end{aligned}
\end{equation}
The analytical stress field is computed from the velocity and pressure fields using
the Cauchy stress relation \eqref{eq:constitutive}. The corresponding prescribed body force is given by
\begin{equation}
\vec f_b(x,y)=
\begin{pmatrix}
\mu\pi^3 \cos(\pi y)  \left(2\sin(\pi x)
+\dfrac{5}{4}\sin(2\pi x) \right)
-0.9
\\
-\mu\pi^3 \sin(\pi y) \left(2\cos(\pi x) +\dfrac{5}{2}\cos(2\pi x) \right)
-0.2
\end{pmatrix}.
\end{equation}

\noindent The test is performed on the unit square domain $\Omega=(0,1)^2$, considering a Newtonian fluid with dynamic viscosity $\mu=1$. The computational domain is discretized using a triangular \ac{fe} mesh composed of $1{,}179{,}648$ elements. The mesh is kept fixed for all simulations, so that variations in the error are driven by the voxel resolution rather than by mesh refinement. To emulate velocity measurements with different spatial resolutions, the analytical
velocity field is subsampled over a sequence of square voxel grids of size
$N_{\mathrm{vox}}\times N_{\mathrm{vox}}$, with
\begin{equation*}
N_{\mathrm{vox}} \in \{4,6,8,12,16,24,32,48,64,96,128,192,256\}.
\end{equation*}

\begin{figure}[htbp]
    \hspace{1.5cm}
    \begin{subfigure}[b]{0.38\textwidth}
        \centering
        \includegraphics[width=\linewidth]{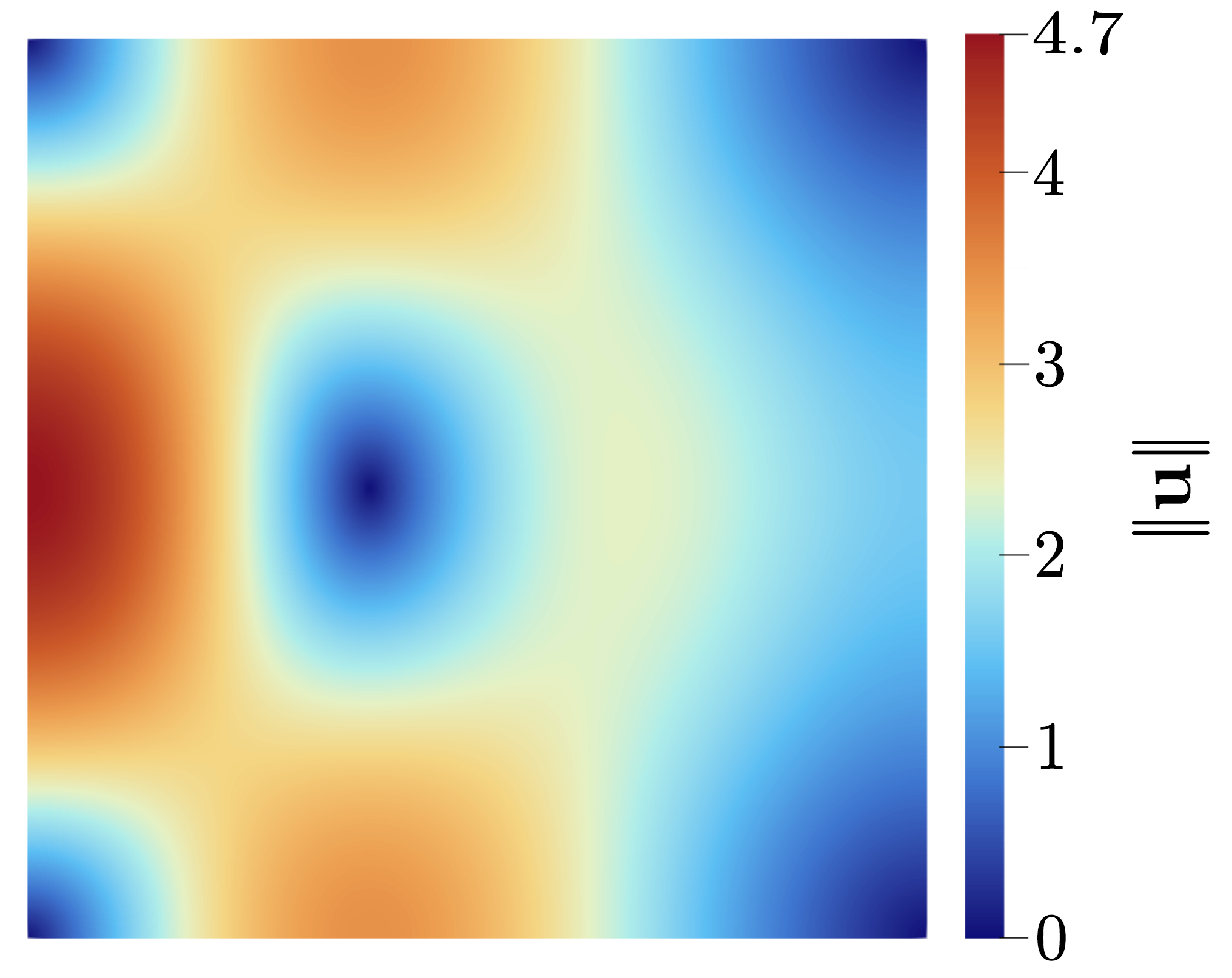}
        \label{fig:velocity_analytic}
    \end{subfigure}
    \hspace{0.5cm}
    \begin{subfigure}[b]{0.4\textwidth}
        \centering
        \includegraphics[width=\linewidth]{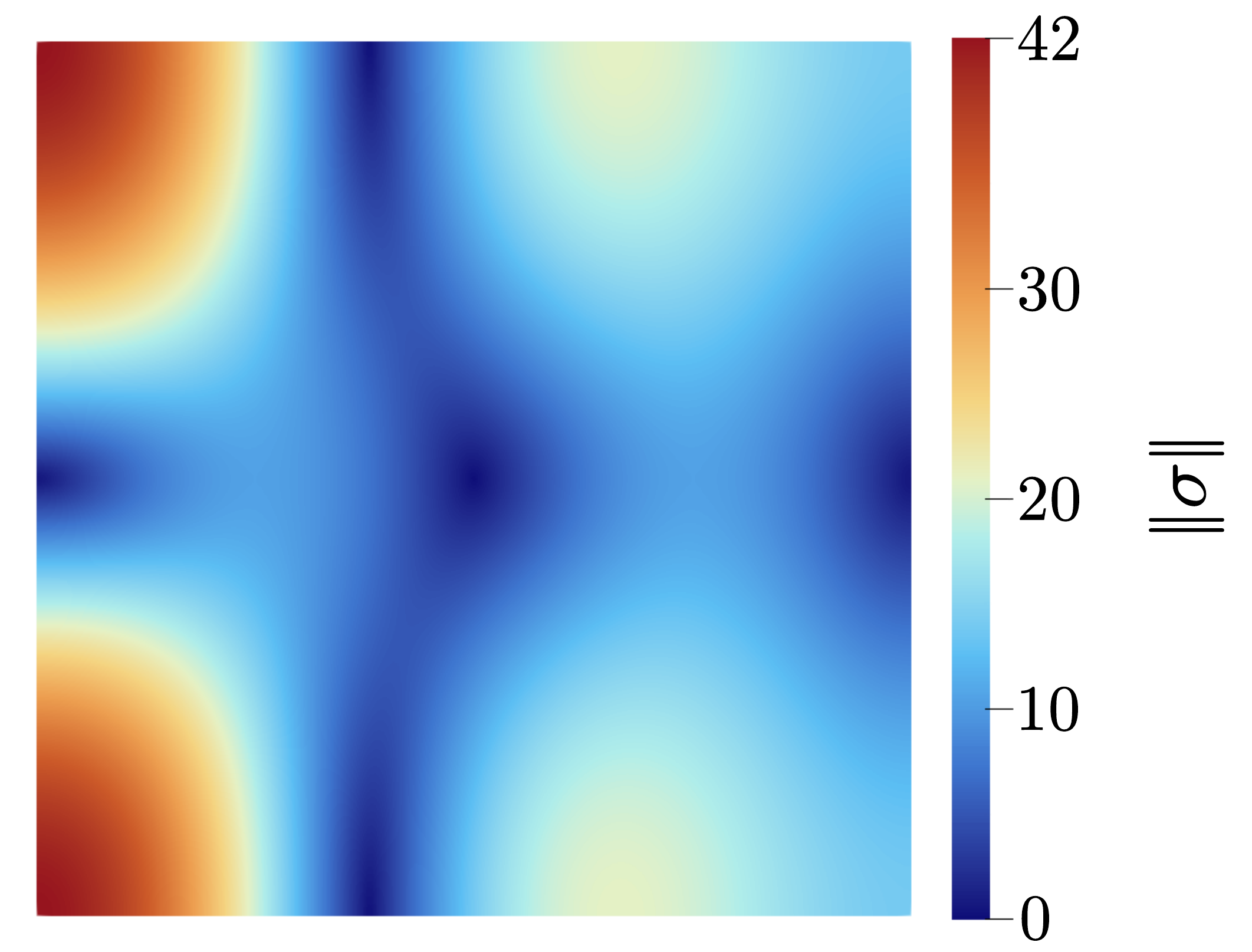}
        \label{fig:stress_analytic}
    \end{subfigure}
    \vspace{-2pt}
    \caption{Analytical velocity and stress fields used in the convergence test.}
    \label{fig:velocity_fields}
\end{figure}

\noindent Therefore, each voxel grid contains $N_{\mathrm{vox}}^2$ voxels and has voxel size $h_{\mathrm{vox}}=1/N_{\mathrm{vox}}$ on the unit square. Following the measurement framework described in Section~\ref{sec:p_recovery}, each voxel defines a local discrete measurement functional based on averaging the available velocity
data inside the voxel. The associated Riesz representers generate the observation space $W_m$, and the resulting measured velocity $\bu_m$ is represented on the fixed computational mesh as a continuous piecewise linear \ac{fe} field. This field is used as input data in the \ac{psp}. The recovered stress field is computed on the same fixed mesh using lowest-order Raviart–Thomas \ac{fe}. 

\begin{figure}[ht]
\centering
\includegraphics[width=\linewidth]{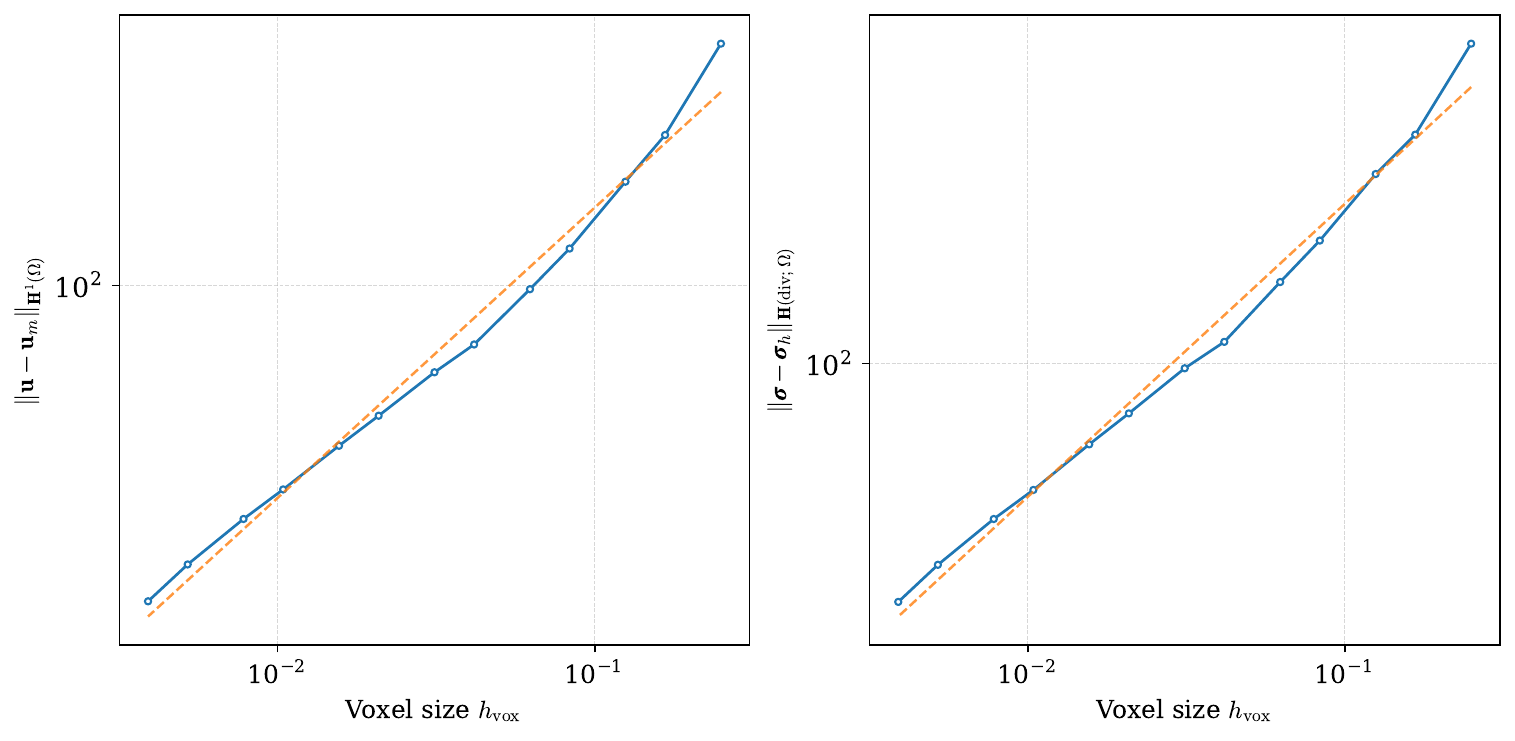}
\caption{Convergence of the velocity measurements and the recovered stress field. Left: $\vec H^1(\Omega)$ error of the velocity field. Right: $\ten H(\mathrm{div};\Omega)$ error of the stress field.}
\label{fig:conv_test_measures}
\end{figure}

\newpage
\noindent We report the errors in the $\vec H^1(\Omega)$ and $\ten H(\dive;\Omega)$ norms.
\begin{equation}
e_{u}(h_{\mathrm{vox}})
=
\|\vec u-\vec u_m\|_{\vec H^1(\Omega)},
\qquad
e_{\sigma}(h_{\mathrm{vox}})
=
\|\ten\sigma-\ten\sigma_h\|_{\ten
H(\dive;\Omega)}.
\end{equation}
The results are shown in Figure~\ref{fig:conv_test_measures}. Both error curves exhibit closely similar convergence behavior, with a least-squares slope close to $0.72$. The agreement between the decay rates indicates that the error in the recovered stress field follows the velocity-measurement error up to a multiplicative constant. This behavior is consistent with the stability estimate \eqref{eq:recovery_error}, derived above for the steady case, which bounds the \ac{psp} reconstruction error in terms of the velocity-measurement error. Therefore, the numerical results
support the theoretical prediction that the \ac{psp} reconstruction is controlled by the quality of the measured velocity field.

\subsection{Venturi channel}
\begin{figure}[ht]
\centering
\includegraphics[width=0.6\linewidth]{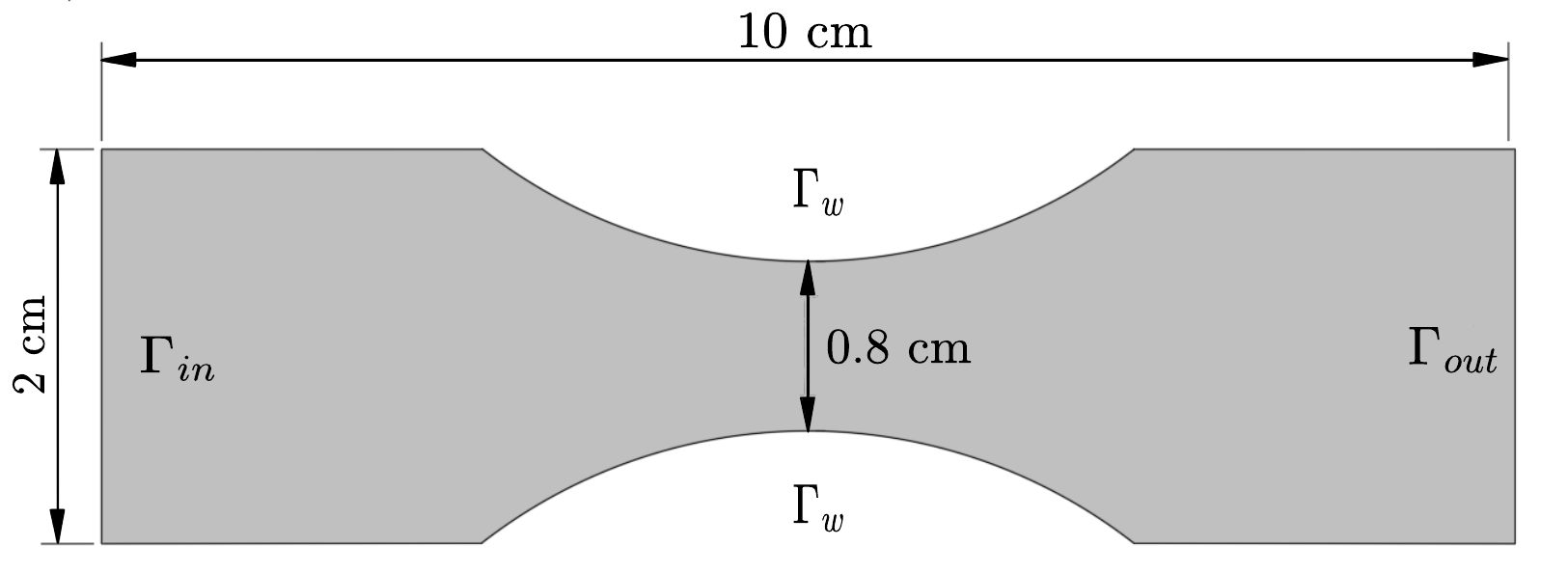}
\caption{Venturi channel geometry}
\label{fig:venturi_domain}
\end{figure}

\begin{figure}[ht]
\centering
\includegraphics[width=0.55\linewidth]{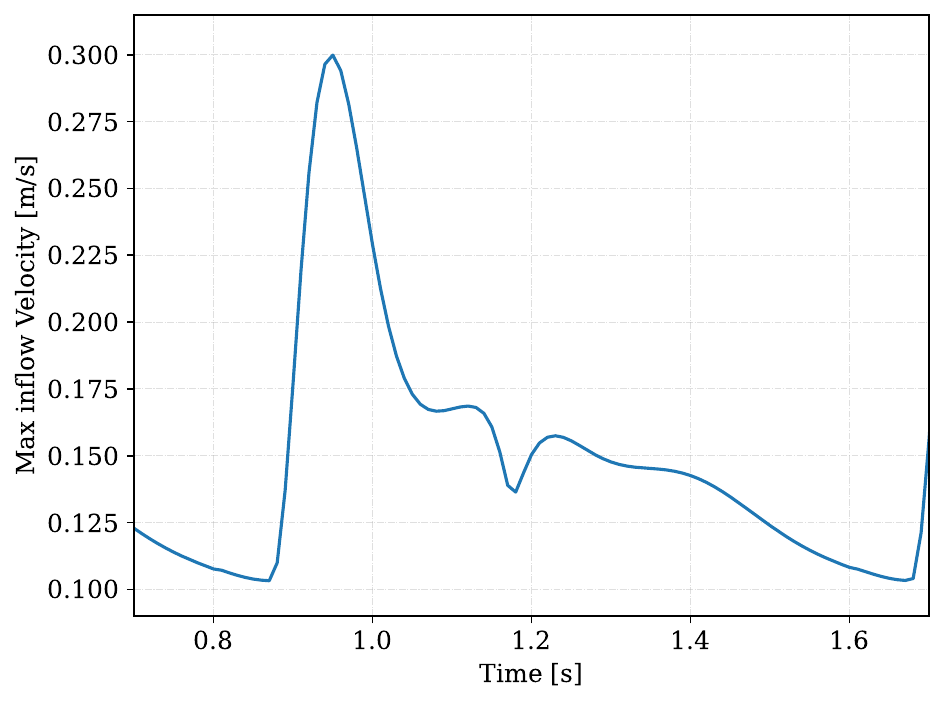}
\caption{Time-dependent maximum inlet velocity $u_{\max}(t)$ per cycle in the Venturi benchmark.}
\label{fig:max_inlet_velocity}
\end{figure}

We next consider a transient Venturi channel benchmark to compare the proposed \ac{psp} with the reference pressure recovery methods described above. The flow is governed by the incompressible \ac{nse} \eqref{eq:NSE} for a Newtonian fluid. The geometry consists of a Venturi channel of length $L=0.1 \ \mathrm{m}$ and inlet
height $H =0.02 \ \mathrm{m}$, as shown in Figure~\ref{fig:venturi_domain}. The boundary is
decomposed as $\partial\Omega=\overline{\Gamma}_{in}\cup\overline{\Gamma}_{out}\cup\overline{\Gamma}_w$.
\noindent The ground-truth velocity and pressure fields are obtained from a direct
numerical simulation in the Venturi domain. The fluid is Newtonian, with dynamic
viscosity $\mu=0.0035 \ \mathrm{Pa}\cdot\mathrm{s}$, density $\rho=1000 \ \mathrm{kg/m^3}$ and no body forces are considered $\bff_b=0$. The reference simulation uses a no-slip condition on the wall $\Gamma_w$, a homogeneous zero traction condition on the outlet $\Gamma_{out}$, and a time-dependent parabolic
velocity profile on the inlet $\Gamma_{in}$
\begin{equation}
    \bu(y,t)
    =
    u_{\text{max}}(t)
    \ \frac{4y(H-y)}{H^2}
    \vec n
    \qquad \text{on } \Gamma_{in} .
\end{equation}
Here, $y$ denotes the transverse coordinate across the inlet section, with
$0\leq y\leq H  $. The function $u_{\max}(t)$ defines a periodic
pulsatile inflow profile used to emulate a cardiac cycle, consistent with the
cardiovascular pressure-recovery applications discussed in the introduction. The
prescribed maximum inlet velocity is shown in Figure~\ref{fig:max_inlet_velocity}.

\noindent With this setup, the ground-truth fields are obtained from a monolithic \ac{fe} simulation of the incompressible \ac{nse} over a time interval of approximately $T=3\,\mathrm{s}$, corresponding to three pulsatile cycles. The spatial
discretization uses stabilized equal-order piecewise linear \ac{fe} for both velocity and pressure, together with a BDF2 time discretization with time step $\Delta t=0.01\,\mathrm{s}$. The nonlinear convective term is treated explicitly through an IMEX strategy, while the viscous and pressure terms are treated implicitly.

\begin{figure}[htbp]
\centering
\begin{subfigure}{0.75\textwidth} 
\centering
\includegraphics[width=\linewidth]{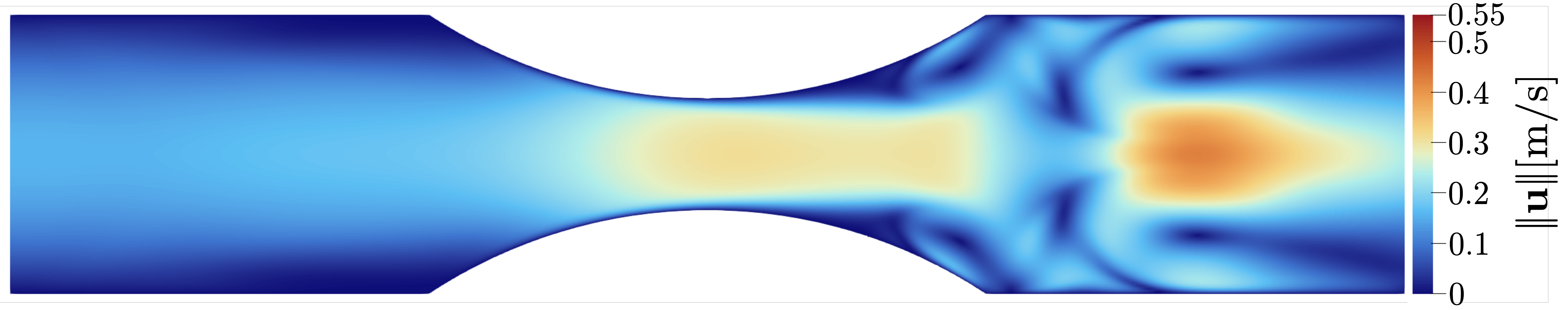}
    \label{fig:gt u field}
\end{subfigure}

\begin{subfigure}{0.75\textwidth}
        \includegraphics[width=\textwidth]{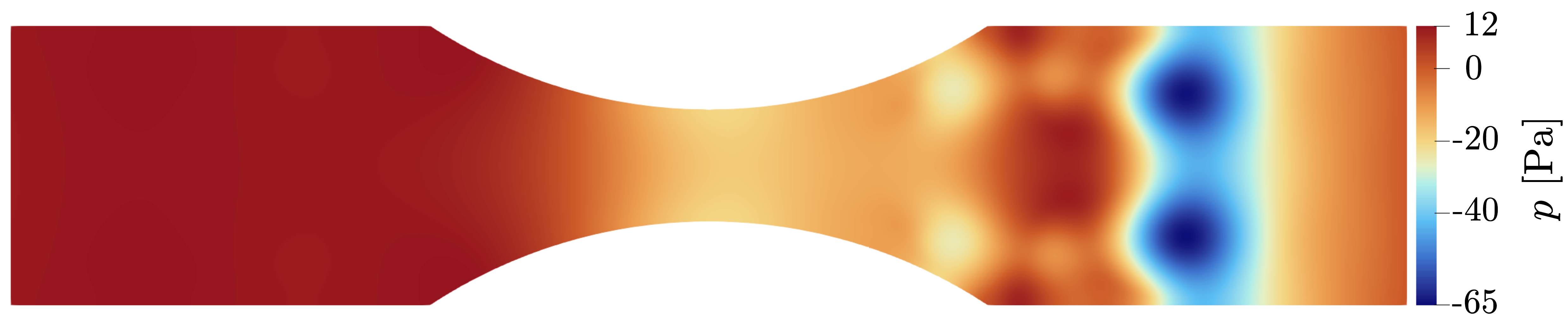}
        \label{fig:gt p field}
    \end{subfigure}

\caption{Venturi channel. Ground truth velocity and pressure field at t=1.25s.}
\label{fig:GT_fields}
\end{figure}

The direct simulation reaches a maximum velocity of approximately $0.55\,\mathrm{m/s}$ and a maximum pressure of about $1.5\,\mathrm{kPa}$. The resulting time-dependent velocity and pressure fields define the ground-truth solution used in the subsequent comparisons, with representative snapshots shown in Figure~\ref{fig:GT_fields}.

\begin{figure}[htbp]
    \centering
    
    \begin{subfigure}{0.75\textwidth}
        \centering
        \includegraphics[width=\textwidth]{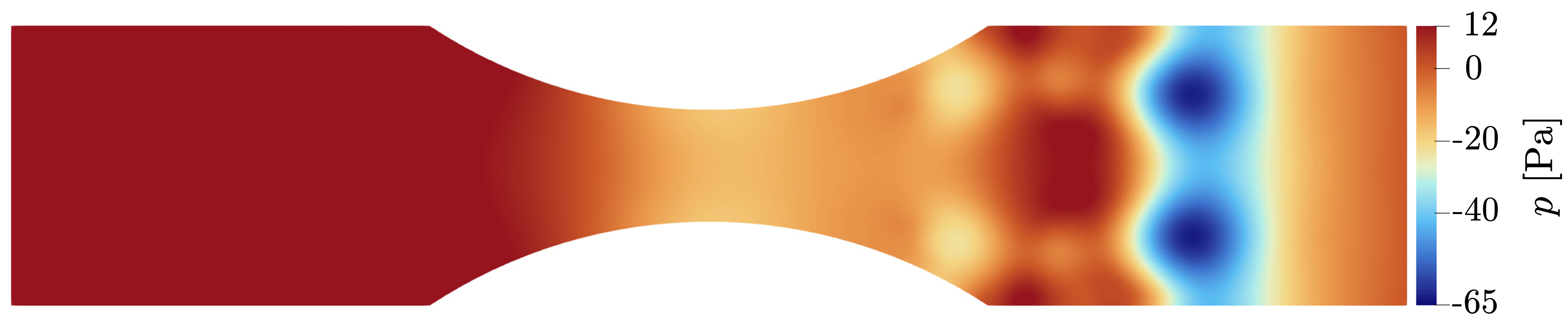}
        \caption{\ac{psp}.}
        \label{fig: ps p field}
    \end{subfigure}

    \begin{subfigure}{0.75\textwidth}
        \centering
        \includegraphics[width=\textwidth]{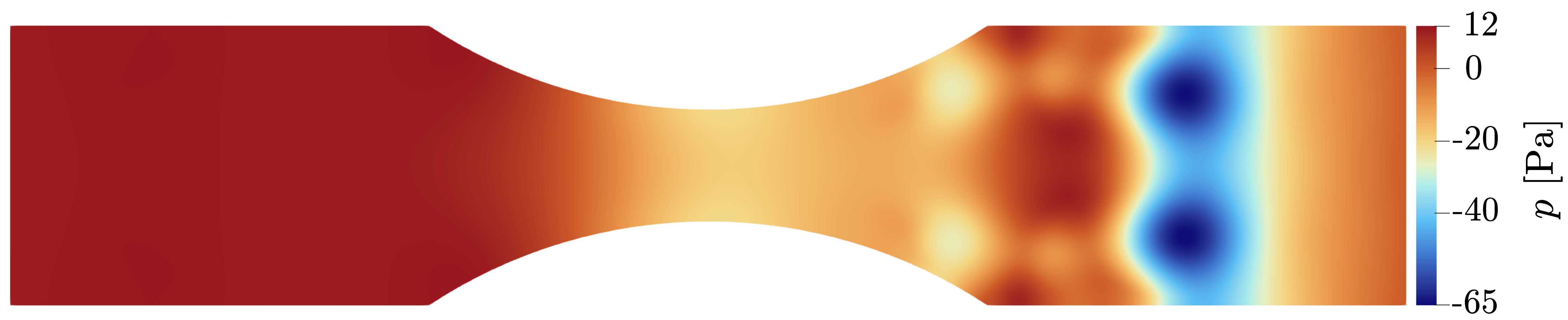}
        \caption{STE.}
        \label{fig: ste p field}
    \end{subfigure}

    \caption{Venturi channel. Recovered pressure fields obtained with the \ac{psp} and \ac{ste} at t=1.25s.}
    \label{fig:gt_ste_pspe_pressure}
\end{figure}

\noindent The ground-truth velocity field is used as input data for all
pressure-recovery estimators, while the ground-truth pressure field provides the
reference solution for evaluating the recovered fields and the corresponding
pressure-drop curves. Figure~\ref{fig:gt_ste_pspe_pressure} shows representative reconstructed pressure fields obtained from this velocity data. The spatial and
temporal setup used to compute these reconstructions is described next.

\newpage
\noindent The recovery problems use the same physical parameters and time step as the reference simulation. The time derivative is approximated using BDF2, consistently with the direct simulation. The stress field is approximated with lowest-order Raviart--Thomas elements on the same Venturi mesh, which contains $106{,}240$ elements. For the \ac{psp} reconstruction, an essential
boundary condition is imposed on the normal stress 
$\ten\sigma\vec n = \left(-p_{\mathrm{ref}}\ten I_d
+ \ten\tau_v(\bu_{\mathrm{GT}}) \right)\vec n$ at the outlet $\Gamma_{out}$. This condition sets a uniform reference pressure $p_{\mathrm{ref}}=0$ at the outlet.

The reconstructed pressure fields are then compared with the ground-truth pressure field. Both reconstructions are visually  indistinguishable from the ground-truth pressure distribution, indicating that the structure of the pressure field is accurately recovered in the noise-free case.

\begin{figure}[htbp]
    \centering
    \includegraphics[width=\textwidth]{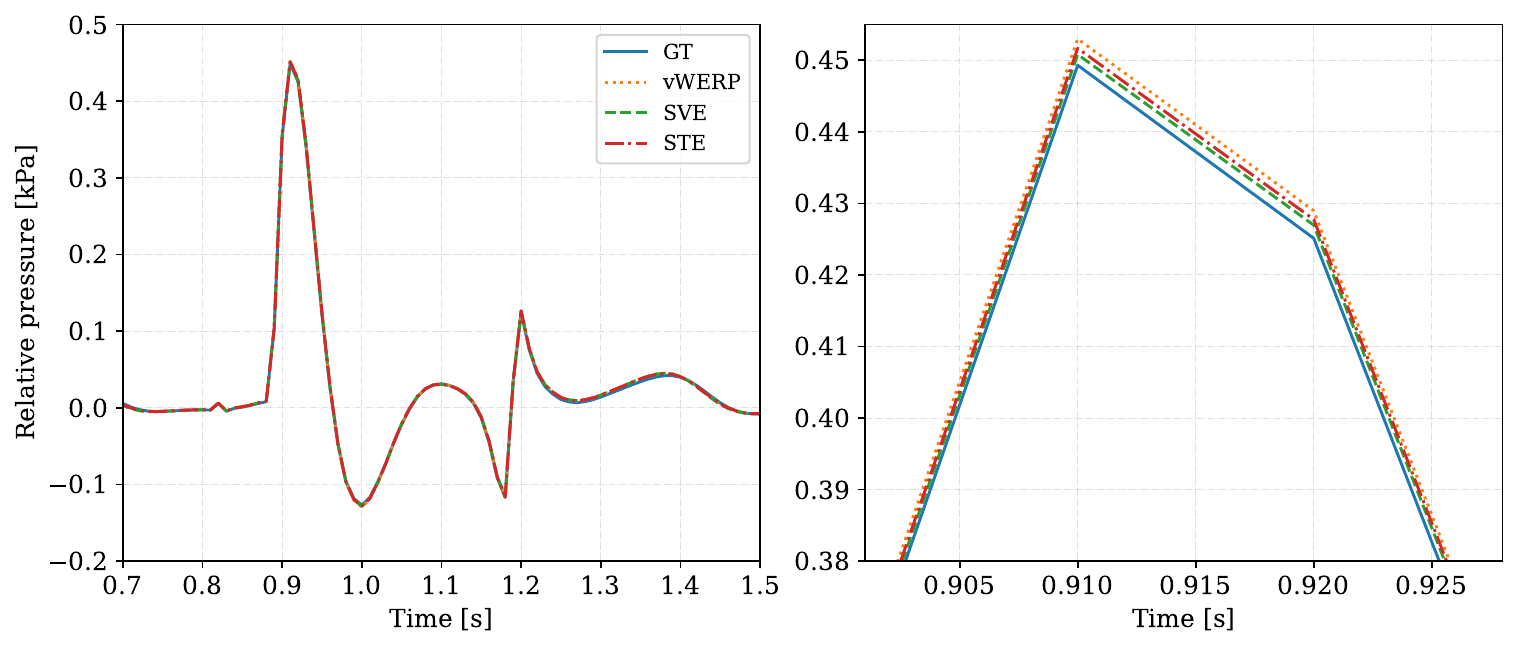}
    \caption{Comparison of relative pressure drop  estimations.}
    \label{fig:pdrop}
\end{figure}

The instantaneous pressure drop, computed according to \eqref{eq:pdrop}, is shown in Figure~\ref{fig:pdrop}. For the field-based estimators, \ac{ste} and \ac{psp}, the pressure drop is computed from the reconstructed pressure fields, whereas \ac{vwerp} directly provides a pressure-drop estimate. All three estimators follow the reference pressure-drop curve throughout the cycle. In this noise-free case, \ac{psp} provides a mild but persistently more accurate pressure-drop prediction than \ac{ste} and \ac{vwerp}.

\subsection{Venturi channel with noisy measurements}
To assess the sensitivity of the pressure-drop estimation methods to
measurement uncertainty, we construct synthetic noisy data from a
ground-truth solution. Let $\boldsymbol{u}(\boldsymbol{x},t)$ denote
the noise-free reference velocity field. The corresponding noisy field is defined as
\begin{equation}
  \boldsymbol{u}^{\text{noisy}}(\boldsymbol{x},t)
  =
  \boldsymbol{u}(\boldsymbol{x},t)
  + \boldsymbol{\eta}(\boldsymbol{x},t),
\end{equation}
where $\boldsymbol{\eta}$ models the measurement noise. The noise is specified using
a standard 4D Flow MRI / PC-MRI parameterization in terms of the velocity encoding
(VENC) and the signal-to-noise ratio (SNR). Let $V_{\text{ENC}}$ be the chosen
velocity encoding value and $\mathrm{SNR}$ the prescribed SNR. The standard deviation of the noise in each velocity component, $\sigma_N$, is given by
\begin{equation}
  \sigma_N
  = \sqrt{\frac{2}{\pi}}\,
    \frac{V_{\text{ENC}}}{\mathrm{SNR}} ,
  \label{eq:sigmaN_venc_snr}
\end{equation}
which corresponds to the propagation of Gaussian magnitude noise to phase-encoded
velocities. In the present work, we use $V_\text{ENC}=0.55\,\mathrm{m/s}$ and
$\mathrm{SNR}=10$, yielding a noise level representative of typical 4D Flow MRI data. For each velocity component, spatial degree of freedom, and time instant, we draw an
independent Gaussian random variable
\begin{equation*}
    \eta \sim \mathcal{N}(0,\sigma_N^2),
\end{equation*}
and assign it to the corresponding component of
$\boldsymbol{\eta}(\boldsymbol{x},t)$. In this way, the components of
$\boldsymbol{\eta}$ are independent and identically distributed with variance
$\sigma_N^2$, and the resulting field $\boldsymbol{u}^{\text{noisy}}$ mimics the
noise level expected for a 4D Flow MRI acquisition with the prescribed
$V_{\text{ENC}}$ and SNR.

\begin{figure}[htbp]
    \centering
    \includegraphics[width=\textwidth]{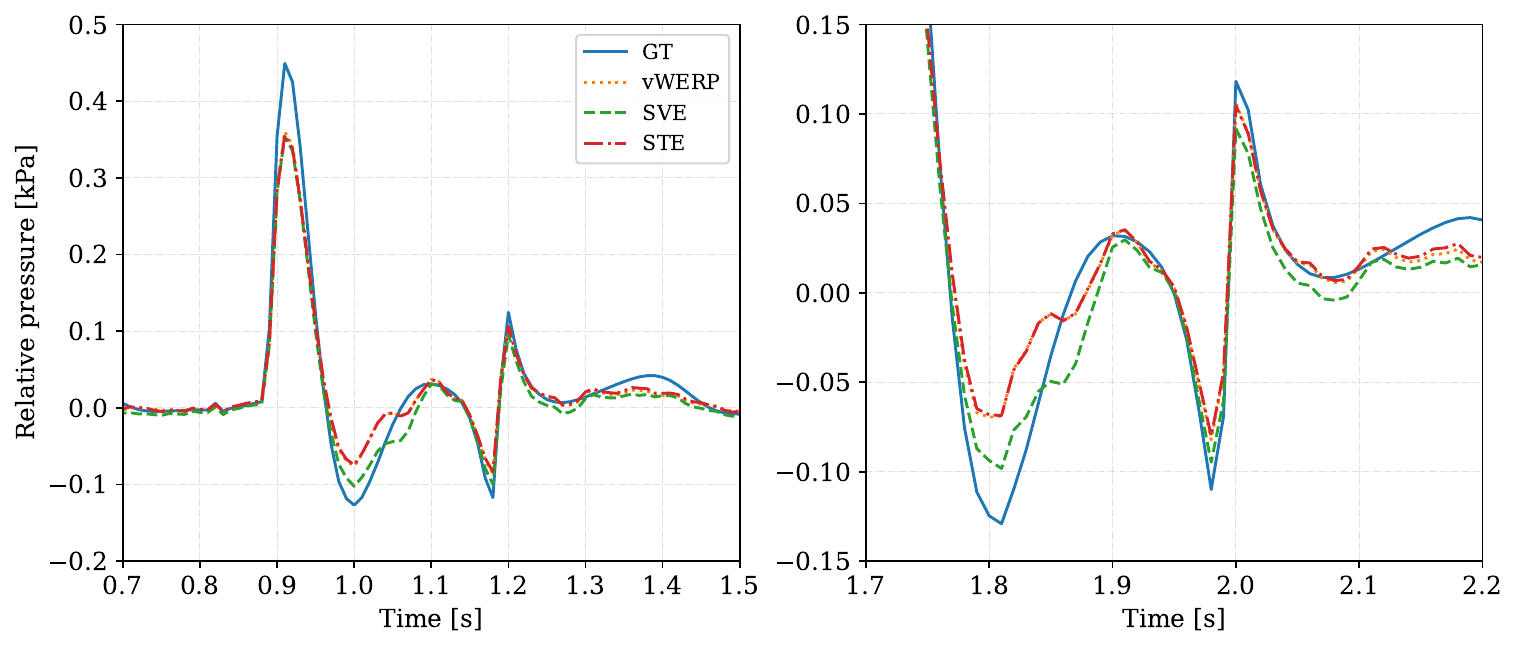}
    \caption{Comparison of relative pressure-drop estimations using noisy velocity measurements.}
    \label{fig:pdrop_noisy_measurements}
\end{figure}

The noisy velocity fields are then used as input data for the three pressure-drop
estimation methods. The numerical setup is the same as in the noise-free Venturi test. The time derivative is approximated using BDF2, whereas the convective and viscous velocity-dependent terms are evaluated using the midpoint approximation
\begin{equation*}
    \boldsymbol{u}^{n+1/2}
    =
    \frac{\boldsymbol{u}^{n+1}+\boldsymbol{u}^{n}}{2}.
\end{equation*}
Thus, $\bu^{n+1/2}$ is used to compute the convective and viscous contributions in
\ac{psp}, \ac{ste}, and \ac{vwerp}. A total of 30 independent noisy simulations are generated and processed with each method. For each simulation, we compute the instantaneous relative pressure drop between the inlet and outlet sections. The final pressure-drop curve reported in Figure \ref{fig:pdrop_noisy_measurements} corresponds to the ensemble-averaged curve over the 30 noisy simulations. The results show that, under noisy velocity measurements, the \ac{psp} estimator provides the closest agreement with the reference pressure-drop curve among the three methods, indicating improved robustness to measurement noise.

\section{Conclusions and perspectives}
\label{sec:conclusions}
In this work, we introduced the \ac{psp}, a stress-based pressure-recovery method. The method reconstructs the full Cauchy stress tensor directly from a known or measured velocity field and then obtains the pressure as a post-processing step from the trace of the recovered stress. In this way, pressure recovery is reformulated as a single-field stress reconstruction problem in $\ten H(\dive;\Omega)$, avoiding auxiliary fields and the saddle-point structure associated with velocity–pressure formulations. Furthermore, the reconstruction holds a divergence-free correction bias that circumvents typical artifacts and spurious measured flow patterns ensuring incompressible mass conservation.

The proposed recovery operator was shown to be well posed under the assumptions
considered in this work. This result provides a stable continuous framework for
stress reconstruction, which is discretized using $\ten H(\dive)$-conforming Raviart--Thomas \ac{fe} applied row-wise to the stress tensor. The numerical convergence test for the \ac{psp} recovery operator showed first-order decay of the
discretization error in both $\ten L^2(\Omega)$ and $\ten H(\dive;\Omega)$ norms, in agreement with the expected behavior of the lowest-order Raviart--Thomas approximation.

The measurement-convergence test further showed that the \ac{psp} estimator can reconstruct stress fields from partial velocity information. The recovered stress error exhibited the same observed convergence rate as the velocity-measurement error, showing that the reconstruction follows the quality of the available velocity data. This behavior is consistent with the stability estimate derived for the measurement-driven problem, which bounds the \ac{psp} reconstruction error in terms of the velocity-measurement error.

For clean Navier--Stokes velocity fields, the transient Venturi benchmark showed that the \ac{psp} accurately recovers the pressure field. The reconstructed pressure distribution was nearly indistinguishable from the ground-truth solution, indicating that the \ac{psp} formulation preserves the spatial structure of the pressure field when the input velocity data are consistent with the governing equations. In terms of relative pressure drop, \ac{psp} provided the closest agreement with the reference curve among the estimators considered.

Under noisy velocity measurements, the \ac{psp} remained robust. In the
noisy Venturi benchmark, the pressure-drop curve obtained with \ac{psp} stayed closer to the ground-truth curve than those obtained with \ac{ste} and \ac{vwerp}. This result suggests improved robustness with respect to measurement noise in the pressure-drop estimation.

Overall, these results indicate that the \ac{psp} formulation provides a stable and accurate framework for pressure recovery from velocity measurements, especially when full pressure fields are required. Future work will extend the proposed estimator to non-Newtonian rheologies and complex three-dimensional geometries with multiple outlets.

\section{Acknowledgments}
FG acknowledges the FONDECYT regular 1261167 FONDECYT regular XXXXX.
\newpage
\bibliographystyle{ieeetr} 
\bibliography{literature}

\end{document}